# SIMULTANEOUS PREDICTION OF INDEPENDENT POISSON OBSERVABLES

By Fumiyasu Komaki

*University of Tokyo*

Simultaneous predictive distributions for independent Poisson observables are investigated. A class of improper prior distributions for Poisson means is introduced. The Bayesian predictive distributions based on priors from the introduced class are shown to be admissible under the Kullback–Leibler loss. A Bayesian predictive distribution based on a prior in this class dominates the Bayesian predictive distribution based on the Jeffreys prior.

**1. Introduction.** Suppose that we have independent observations $x(1), x(2), \dots, x(n)$, where $x(l) = (x_1(l), x_2(l), \dots, x_d(l))$ $(l \in \{1, 2, \dots, n\})$ is a set of $d$ independent Poisson random variables with unknown mean parameters $\lambda_1, \lambda_2, \dots, \lambda_d$. We write $x^{(n)} = (x(1), x(2), \dots, x(n))$ and $\lambda = (\lambda_1, \lambda_2, \dots, \lambda_d)$. An unobserved set $x_{(m)} = (x(n+1), x(n+2), \dots, x(n+m))$ from the same distribution is predicted by using a predictive distribution $\hat{p}(x_{(m)}; x^{(n)})$. We adopt the Kullback–Leibler divergence from the true distribution to a predictive distribution,

$$D(p(x_{(m)}|\lambda), \hat{p}(x_{(m)}; x^{(n)})) = \sum_{x_{(m)}} p(x_{(m)}|\lambda) \log \frac{p(x_{(m)}|\lambda)}{\hat{p}(x_{(m)}; x^{(n)})}, \tag{1}$$

which has a natural information theoretic meaning, as a loss function.

By sufficiency reduction, it suffices to consider the problem of predicting $y = (y_1, y_2, \dots, y_d)$ using $x = (x_1, x_2, \dots, x_d)$, where

$$x = \sum_{i=1}^{n} x(i) = \left( \sum_{i=1}^{n} x_1(i), \sum_{i=1}^{n} x_2(i), \dots, \sum_{i=1}^{n} x_d(i) \right),$$

$$y=\sum_{j=1}^{m}x(n+j)=\left(\sum_{j=1}^{m}x_1(n+j),\sum_{j=1}^{m}x_2(n+j),\dots,\sum_{j=1}^{m}x_d(n+j)\right),$$

under the loss

$$D(p(y|\lambda),\hat{p}(y;x))=\sum_{y}p(y|\lambda)\log\frac{p(y|\lambda)}{\hat{p}(y;x)}. \tag{2}$$

In the following, we assume that $x=(x_1,x_2,\dots,x_d)$ and $y=(y_1,y_2,\dots,y_d)$ are distributed according to

$$\begin{aligned}p(x|\lambda)&=\prod_{i=1}^{d}p(x_i|\lambda)\\&=\exp\{-(a\lambda_1+a\lambda_2+\cdots+a\lambda_d)\}\frac{(a\lambda_1)^{x_1}}{x_1!}\frac{(a\lambda_2)^{x_2}}{x_2!}\cdots\frac{(a\lambda_d)^{x_d}}{x_d!}\end{aligned}$$

and

$$\begin{aligned}p(y|\lambda)&=\prod_{i=1}^{d}p(y_i|\lambda)\\&=\exp\{-(b\lambda_1+b\lambda_2+\cdots+b\lambda_d)\}\frac{(b\lambda_1)^{y_1}}{y_1!}\frac{(b\lambda_2)^{y_2}}{y_2!}\cdots\frac{(b\lambda_d)^{y_d}}{y_d!},\end{aligned}$$

respectively, and that $a$ and $b$ are positive real numbers.

There exist many studies that recommend using Bayesian predictive densities of the form

$$p_\pi(x_{(m)}|x^{(n)})=\frac{\int p(x_{(m)}|\theta)p(x^{(n)}|\theta)\pi(\theta)\,d\theta}{\int p(x^{(n)}|\bar{\theta})\pi(\bar{\theta})\,d\bar{\theta}},$$

rather than plug-in densities of the form $p(x_{(m)}|\hat{\theta})$, where $\{p(x|\theta)|\theta\in\Theta\}$ is a parametric model, $\pi(\theta)$ is a prior and $\hat{\theta}$ is an estimate of $\theta$; see Aitchison and Dunsmore (1975) and Geisser (1993).

When we use a Bayesian procedure, the choice of a prior distribution is an important problem. Noninformative prior distributions or vague prior distributions are often used to construct Bayesian predictive distributions. The Jeffreys prior naturally arises from various discussions based on the Kullback–Leibler divergence [see Hartigan (1965), Akaike (1978), Bernardo (1979) and Clarke and Barron (1994)]. However, Bayesian methods based on the Jeffreys prior do not always perform satisfactorily, especially in problems with multidimensional parameters [see, e.g., Jeffreys (1961), page 182, and Berger and Bernardo (1989)].

Here, we investigate the use of shrinkage priors, which give more weight to parameter values close to the origin than the Jeffreys prior does, for constructing predictive distributions dominating the predictive distribution based on the Jeffreys prior. If we adopt a plug-in distribution $p(y|\hat{\lambda}(x))$ as a predictive distribution, the loss (1) for the plug-in distribution can be regarded as a loss for the estimator $\hat{\theta}$. Thus, predictive distribution theory is a natural generalization of estimation theory under the Kullback–Leibler loss.

Since Stein (1956) showed that the maximum likelihood estimator for the mean vector of the $d$-dimensional Normal model $N_d(\mu, I)$ is not admissible when $d \geq 3$ and James and Stein (1961) introduced an estimator dominating the maximum likelihood estimator, numerous studies have been done on shrinkage methods for parameter estimation.

For the means of $d$ independent Poisson distributions, Clevenson and Zidek (1975) proposed a class of estimators dominating the maximum likelihood estimator when $d \geq 2$ under the normalized square loss $\sum_i (\hat{\lambda}_i - \lambda_i)^2/\lambda_i$. Many studies on simultaneous estimation of Poisson means have been done under the loss function $\sum_i (\hat{\lambda}_i - \lambda_i)^2/{\lambda_i}^m$, where $m$ is a nonnegative integer.

Ghosh and Yang (1988) characterized linear admissible estimators of the form $\hat{\lambda}_i = c_i x_i + b_i$ under the Kullback–Leibler loss $D(p(y|\lambda), p(y|\hat{\lambda}(x)))$. There are relatively few studies of estimation under the Kullback–Leibler loss compared with the number of studies based on other loss functions such as squared-error. What is called Stein's loss is the Kullback–Leibler divergence with the direction opposite to our setting (1).

In contrast to the large number of studies on parameter estimation, little attention has been given to decision theory of predictive distributions except for some studies on group models [Murray (1977) and Ng (1980)] and some recent work from an asymptotic viewpoint [Vidoni (1995), Komaki (1996) and Haussler and Opper (1997)]. In particular, it seems that no studies have been done on the admissibility of predictive distributions. Recently, however, Komaki (2001) considered the $d$-dimensional Normal model $N_d(\mu, I)$, $d \geq 3$, and showed that the Bayesian predictive distribution based on Stein's harmonic prior $\pi_{\mathrm{S}}(\mu) \propto \|\mu\|^{-(d-2)}$ [Stein (1974)] incorporates the advantage of shrinkage methods and dominates the Bayesian predictive distribution based on the Lebesgue prior $\pi_{\mathrm{I}}(\mu) \propto 1$, which is the best predictive distribution invariant under the translation group. Since a lot of statistical problems are naturally formulated as prediction problems [Aitchison and Dunsmore (1975) and Geisser (1993)], this kind of approach seems to be useful for many problems, and further decision theoretic studies especially on admissibility are required.

In Section 2 we introduce a class of improper prior densities for Poisson means and show that the predictive distributions based on the proposed priors are admissible under the Kullback–Leibler loss. In Section 3 we show that a Bayesian predictive distribution based on a prior $\pi_{\mathrm{S}}(\lambda)$ in the introduced class dominates the Bayesian predictive distribution based on the Jeffreys prior, and that the plug-in distribution with the generalized Bayes estimator based on $\pi_{\mathrm{S}}(\lambda)$ is inadmissible under the Kullback–Leibler loss. In Section 4 we discuss the relation between the main results here and several previous studies on Bayesian theory from asymptotic viewpoints.

**2. A class of admissible predictive distributions.** We introduce a class of improper prior densities,

$$\pi_{\alpha,\beta}(\lambda)\,d\lambda_1\,d\lambda_2\cdots d\lambda_d \propto \frac{\lambda_1^{\beta_1-1}\lambda_2^{\beta_2-1}\cdots\lambda_d^{\beta_d-1}}{(\lambda_1+\lambda_2+\cdots+\lambda_d)^{\alpha}}\,d\lambda_1\,d\lambda_2\cdots d\lambda_d \tag{3}$$

with $0 < -\alpha + \sum_i \beta_i \le 1$ and $\beta_i > 0$, $i = 1, 2, \ldots, d$.

THEOREM 1. *The Bayesian predictive distribution based on the prior*

$$\pi_{\alpha,\beta}(\lambda)\,d\lambda_1\,d\lambda_2\cdots d\lambda_d \propto \frac{\lambda_1^{\beta_1-1}\lambda_2^{\beta_2-1}\cdots\lambda_d^{\beta_d-1}}{(\lambda_1+\lambda_2+\cdots+\lambda_d)^{\alpha}}\,d\lambda_1\,d\lambda_2\cdots d\lambda_d$$

*with* $-\alpha + \sum_i \beta_i > 0$ *and* $\beta_i > 0$, $i = 1, 2, \ldots, d$, *is given by*

$$\begin{aligned} &p_{\pi_{\alpha,\beta}}(y|x)\\ &\quad= \left(\frac{a}{a+b}\right)^{\sum_i x_i - \alpha + \sum_i \beta_i}\left(\frac{b}{a+b}\right)^{\sum_i y_i}\\ &\qquad\times \frac{\Gamma(\sum_i x_i + \sum_i y_i - \alpha + \sum_i \beta_i)\Gamma(\sum_i x_i + \sum_i \beta_i)}{\Gamma(\sum_i x_i - \alpha + \sum_i \beta_i)\Gamma(\sum_i x_i + \sum_i y_i + \sum_i \beta_i)}\\ &\qquad\times \frac{\Gamma(x_1+y_1+\beta_1)\Gamma(x_2+y_2+\beta_2)\cdots\Gamma(x_d+y_d+\beta_d)}{\Gamma(x_1+\beta_1)\Gamma(x_2+\beta_2)\cdots\Gamma(x_d+\beta_d)y_1!y_2!\cdots y_d!}. \end{aligned}$$

PROOF. By using Lemma 1 below, we have

$$\begin{aligned} &p_{\pi_{\alpha,\beta}}(y|x)\\ &\quad= \frac{\int \pi_{\alpha,\beta}(\lambda)\prod_{i=1}^d\{\exp(-a\lambda_i)(a\lambda_i)^{x_i}/x_i!\}\prod_{j=1}^d\{\exp(-b\lambda_j)(b\lambda_j)^{y_j}/y_j!\}\,d\lambda}{\int \pi_{\alpha,\beta}(\bar\lambda)\prod_{k=1}^d\{\exp(-a\bar\lambda_k)(a\bar\lambda_k)^{x_k}/x_k!\}\,d\bar\lambda}\\ &\quad= \frac{\int \pi_{\alpha,\beta}(\lambda)\prod_{i=1}^d[\exp\{-(a+b)\lambda_i\}\{(a+b)\lambda_i\}^{x_i+y_i}]\,d\lambda}{\int \pi_{\alpha,\beta}(\bar\lambda)\prod_{k=1}^d\{\exp(-a\bar\lambda_k)(a\bar\lambda_k)^{x_k}\}\,d\bar\lambda}\\ &\qquad\times\prod_{j=1}^d\frac{a^{x_j}b^{y_j}}{(a+b)^{x_j+y_j}y_j!}\\ &\quad= \frac{a^{\sum_i x_i - \alpha + \sum_i \beta_i}b^{\sum_i y_i}}{(a+b)^{\sum_j x_j + \sum_j y_j - \alpha + \sum_j \beta_j}} \end{aligned}$$

$$
\begin{aligned}
&\times \frac{\Gamma(\sum_i x_i + \sum_i y_i - \alpha + \sum_i \beta_i)}{\Gamma(\sum_i x_i + \sum_i y_i + \sum_i \beta_i)} \\
&\times \Gamma(x_1 + y_1 + \beta_1)\Gamma(x_2 + y_2 + \beta_2) \cdots \Gamma(x_d + y_d + \beta_d) \\
&\times \Bigg[ \frac{\Gamma(\sum_i x_i - \alpha + \sum_i \beta_i)}{\Gamma(\sum_i x_i + \sum_i \beta_i)} \\
&\qquad \times \Gamma(x_1 + \beta_1)\Gamma(x_2 + \beta_2) \cdots \Gamma(x_d + \beta_d) y_1! y_2! \cdots y_d! \Bigg]^{-1}.
\end{aligned}
$$

Thus we obtain the desired result. □

LEMMA 1. *When $-\alpha + \sum_i \beta_i > 0$ and $\beta_i > 0$, $i = 1, 2, \ldots, d$, we have that*

$$
\begin{aligned}
&\int \frac{\lambda_1^{\beta_1 - 1} \lambda_2^{\beta_2 - 1} \cdots \lambda_d^{\beta_d - 1}}{(\lambda_1 + \lambda_2 + \cdots + \lambda_d)^\alpha} \prod_{i=1}^d \{\exp(-a\lambda_i)(a\lambda_i)^{x_i}\} \, d\lambda_1 \, d\lambda_2 \cdots d\lambda_d \\
&\qquad = a^{\alpha - \sum_i \beta_i} \frac{\Gamma(\sum_i x_i - \alpha + \sum_i \beta_i)}{\Gamma(\sum_i x_i + \sum_i \beta_i)} \prod_{i=1}^d \Gamma(x_i + \beta_i).
\end{aligned}
$$

The proof of Lemma 1 is given in the Appendix.

Let $\mathcal{P}$ be the class of predictive distributions that have finite risk for all values of $\lambda$. For example, the plug-in distribution

$$
p(y|\hat{\lambda}(x)) = \prod_{i=1}^d \exp\left(-\frac{bx_i}{a}\right) \frac{(bx_i/a)^{y_i}}{y_i!}
$$

with the maximum likelihood estimator $\hat{\lambda}(x) = x/a$ is not included in $\mathcal{P}$, because the loss (2) becomes infinite when $\lambda_i > 0$ and $\hat{\lambda}_i(x) = 0$. If a predictive distribution is admissible in $\mathcal{P}$, then it is admissible in the class of all predictive distributions.

Before proving the admissibility of the proposed class of Bayesian predictive distributions, we establish the following theorem.

THEOREM 2. *If $\hat{p}(y; x) \in \mathcal{P}$, then the risk function*

$$
r_{\hat{p}}(\lambda) = E[D(p(y|\lambda), \hat{p}(y; x))|\lambda]
$$

*is a continuous function of $\lambda$.*

PROOF. The risk function is given by

$$
\sum_x p(x|\lambda) \sum_y p(y|\lambda) \log \frac{p(y|\lambda)}{\hat{p}(y; x)}
$$

$$= \sum_y p(y|\lambda) \log p(y|\lambda) \tag{4}$$

$$+ \sum_x p(x|\lambda) \sum_y p(y|\lambda)\{-\log \hat{p}(y;x)\}.$$

The first term on the right-hand side of (4) is

$$\sum_y p(y|\lambda) \log p(y|\lambda)$$

$$= \sum_{i=1}^d \left[ \sum_{y_i=0}^\infty \exp(-b\lambda_i) \frac{(b\lambda_i)^{y_i}}{y_i!} \{-b\lambda_i + y_i \log(b\lambda_i) - \log y_i!\} \right].$$

This is finite for all values of $\lambda$ and is a continuous function of $\lambda$. The second term on the right-hand side of (4) is

$$\sum_x p(x|\lambda) \sum_y p(y|\lambda)\{-\log \hat{p}(y;x)\}$$

$$= \sum_x \sum_y \prod_{i=1}^d \exp(-a\lambda_i) \frac{(a\lambda_i)^{x_i}}{x_i!}$$

$$\times \prod_{j=1}^d \exp(-b\lambda_j) \frac{(b\lambda_j)^{y_j}}{y_j!} \{-\log \hat{p}(y;x)\}$$

$$= \exp\left\{-(a+b) \sum_i \lambda_i\right\} \tag{5}$$

$$\times \left[ \sum_x \sum_y \frac{a^{\sum x_i} b^{\sum y_j}}{x_1! x_2! \cdots x_d! y_1! y_2! \cdots y_d!} \{-\log \hat{p}(y;x)\} \right.$$

$$\left. \times \lambda_1^{x_1+y_1} \lambda_2^{x_2+y_2} \cdots \lambda_d^{x_d+y_d} \right].$$

If $\hat{p}(y;x) \in \mathcal{P}$, the power series in $\lambda_1, \lambda_2, \ldots, \lambda_d$,

$$\sum_x \sum_y \frac{a^{\sum x_i} b^{\sum y_j}}{x_1! x_2! \cdots x_d! y_1! y_2! \cdots y_d!} \{-\log \hat{p}(y;x)\} \lambda_1^{x_1+y_1} \lambda_2^{x_2+y_2} \cdots \lambda_d^{x_d+y_d},$$

converges absolutely for all $\lambda \in \mathbf{R}^d$. Thus, (5) is a continuous function of $\lambda$. Therefore, the risk function is continuous for all values of $\lambda$ if $\hat{p}(y;x) \in \mathcal{P}$. $\square$

THEOREM 3. *For every $d \geq 1$, the Bayesian predictive distributions based on the priors in the class $\{\pi_{\alpha,\beta}(\lambda) : 0 < -\alpha + \sum_{i=1}^d \beta_i \leq 1,\ \beta_i > 0,\ i = 1, 2, \ldots, d\}$ defined by (3) are admissible under the Kullback–Leibler loss.*

The proof of Theorem 3 is given in the Appendix.

## 3. A shrinkage prior dominating the Jeffreys prior.

In this section, we show that the Bayesian predictive distribution based on the Jeffreys prior is inadmissible and give an explicit form of a shrinkage predictive distribution dominating the Bayesian predictive distribution based on the Jeffreys prior.

First, we show that the Bayesian predictive distribution $p_{\pi_{\alpha,\beta}}(y|x)$ is inadmissible when $-\alpha + \sum_i \beta_i > 1$.

THEOREM 4. *When $-\alpha+\sum_i \beta_i > 1$ and $\beta_i > 0$, $i = 1, 2, \ldots, d$, the Bayesian predictive distribution $p_{\pi_{\alpha,\beta}}(y|x)$ based on $\pi_{\alpha,\beta}(\lambda)$ is dominated by the Bayesian predictive distribution $p_{\pi_{\tilde{\alpha},\tilde{\beta}}}(y|x)$ based on $\pi_{\tilde{\alpha},\tilde{\beta}}(\lambda)$, where $\tilde{\alpha} := \sum_i \beta_i - 1$ and $\tilde{\beta} = (\tilde{\beta}_1, \tilde{\beta}_2, \ldots, \tilde{\beta}_d) := (\beta_1, \beta_2, \ldots, \beta_d)$.*

PROOF. From Theorem 1, we have

$$
\begin{aligned}
E[D(&p(y|\lambda), p_{\pi_{\alpha,\beta}}(y|x))|\lambda] - E[D(p(y|\lambda), p_{\pi_{\tilde{\alpha},\tilde{\beta}}}(y|x))|\lambda] \\
&= E\left[\log \frac{p_{\pi_{\tilde{\alpha},\tilde{\beta}}}(y|x)}{p_{\pi_{\alpha,\beta}}(y|x)}\middle|\lambda\right] \\
&= E\left[\log\left\{\left(\frac{a}{a+b}\right)^{-\tilde{\alpha}+\alpha}\right.\right. \\
&\qquad\qquad \left.\left.\times \frac{\Gamma(\sum x_i + \sum y_i - \tilde{\alpha} + \sum \beta_i)\Gamma(\sum x_i - \alpha + \sum \beta_i)}{\Gamma(\sum x_i - \tilde{\alpha} + \sum \beta_i)\Gamma(\sum x_i + \sum y_i - \alpha + \sum \beta_i)}\right\}\middle|\lambda\right] \\
&= E\left[\log\Gamma\left(\sum_i x_i - \alpha + \sum_i \beta_i\right) - \log\Gamma\left(\sum_i x_i - \tilde{\alpha} + \sum_i \beta_i\right)\right. \\
&\qquad\qquad\qquad\qquad \left. - (\tilde{\alpha} - \alpha)\log a \middle|\lambda\right] \\
&\qquad - E\left[\log\Gamma\left(\sum_i x_i + \sum_i y_i - \alpha + \sum_i \beta_i\right)\right. \\
&\qquad\qquad - \log\Gamma\left(\sum_i x_i + \sum_i y_i - \tilde{\alpha} + \sum_i \beta_i\right) \\
&\qquad\qquad\qquad \left. - (\tilde{\alpha} - \alpha)\log(a+b)\middle|\lambda\right] \qquad (6) \\
&= E\left[\log\Gamma\left(\sum_i x_i + 1 + \tilde{\alpha} - \alpha\right) - \log\Gamma\left(\sum_i x_i + 1\right)\right.
\end{aligned}
$$

$$- (\tilde{\alpha} - \alpha) \log a \bigg| \lambda \bigg]$$

$$- E\left[\log \Gamma\left(\sum_i x_i + \sum_i y_i + 1 + \tilde{\alpha} - \alpha\right) - \log \Gamma\left(\sum_i x_i + \sum_i y_i + 1\right) - (\tilde{\alpha} - \alpha) \log(a+b) \bigg| \lambda \right].$$

When $\mu := \sum_i \lambda_i = 0$,

$$\begin{aligned} &E[D(p(y|\lambda), p_{\pi_{\alpha,\beta}}(y|x))|\lambda] - E[D(p(y|\lambda), p_{\pi_{\tilde{\alpha},\tilde{\beta}}}(y|x))|\lambda] \\ &\quad = (\tilde{\alpha} - \alpha) \log \frac{a+b}{a} > 0. \end{aligned} \tag{7}$$

When $\mu > 0$, by using Lemma 2 below, we have

$$E[D(p(y|\lambda), p_{\pi_{\alpha,\beta}}(y|x))|\lambda] - E[D(p(y|\lambda), p_{\pi_{\tilde{\alpha},\tilde{\beta}}}(y|x))|\lambda] > 0,$$

since $\sum_i x_i + \sum_i y_i$ and $\sum_i x_i$ are Poisson random variables with parameters $(a+b) \sum_i \lambda_i$ and $a \sum_i \lambda_i$, respectively. □

LEMMA 2. *Let $X$ be a Poisson random variable with mean $\mu$. Then*

$$E[\log \Gamma(X + 1 + c) - \log \Gamma(X + 1) - c \log \mu | \mu],$$

*where $c$ is a positive constant, is a strictly decreasing function of $\mu > 0$.*

The proof of Lemma 2 is given in the Appendix.

In the following, we set $\pi_{\mathrm{S}}(\lambda) = \pi_{\alpha=d/2-1,\beta=(1/2,\dots,1/2)}(\lambda)$.

COROLLARY 1. *When $d \geq 3$, the Bayesian predictive distribution $p_{\pi_{\mathrm{S}}}(y|x)$ based on the prior $\pi_{\mathrm{S}}(\lambda)$ dominates the Bayesian predictive distribution $p_{\pi_{\mathrm{J}}}(y|x)$ based on the Jeffreys prior*

$$\pi_{\mathrm{J}}(\lambda)\, d\lambda_1\, d\lambda_2 \cdots d\lambda_d \propto \frac{1}{(\lambda_1 \lambda_2 \cdots \lambda_d)^{1/2}}\, d\lambda_1\, d\lambda_2 \cdots d\lambda_d.$$

PROOF. The Jeffreys prior is equal to $\pi_{\alpha=0,\beta=(1/2,\dots,1/2)}$. The desired results follow from Theorem 4 because $-\alpha + \sum_i \beta_i = d/2 > 1$. □

Figure 1 shows the difference between the risk of $p_{\pi_{\mathrm{J}}}(y|x)$ and that of $p_{\pi_{\mathrm{S}}}(y|x)$.

Since $p_{\pi_{\mathrm{S}}}(y|x)$, based on the prior $\pi_{\mathrm{S}}(\lambda)$, dominates $p_{\pi_{\mathrm{J}}}(y|x)$, based on the Jeffreys prior $\pi_{\mathrm{J}}(\lambda)$, it seems to be reasonable to adopt $\pi_{\mathrm{S}}(\lambda)$ as a default prior instead of $\pi_{\mathrm{J}}(\lambda)$.

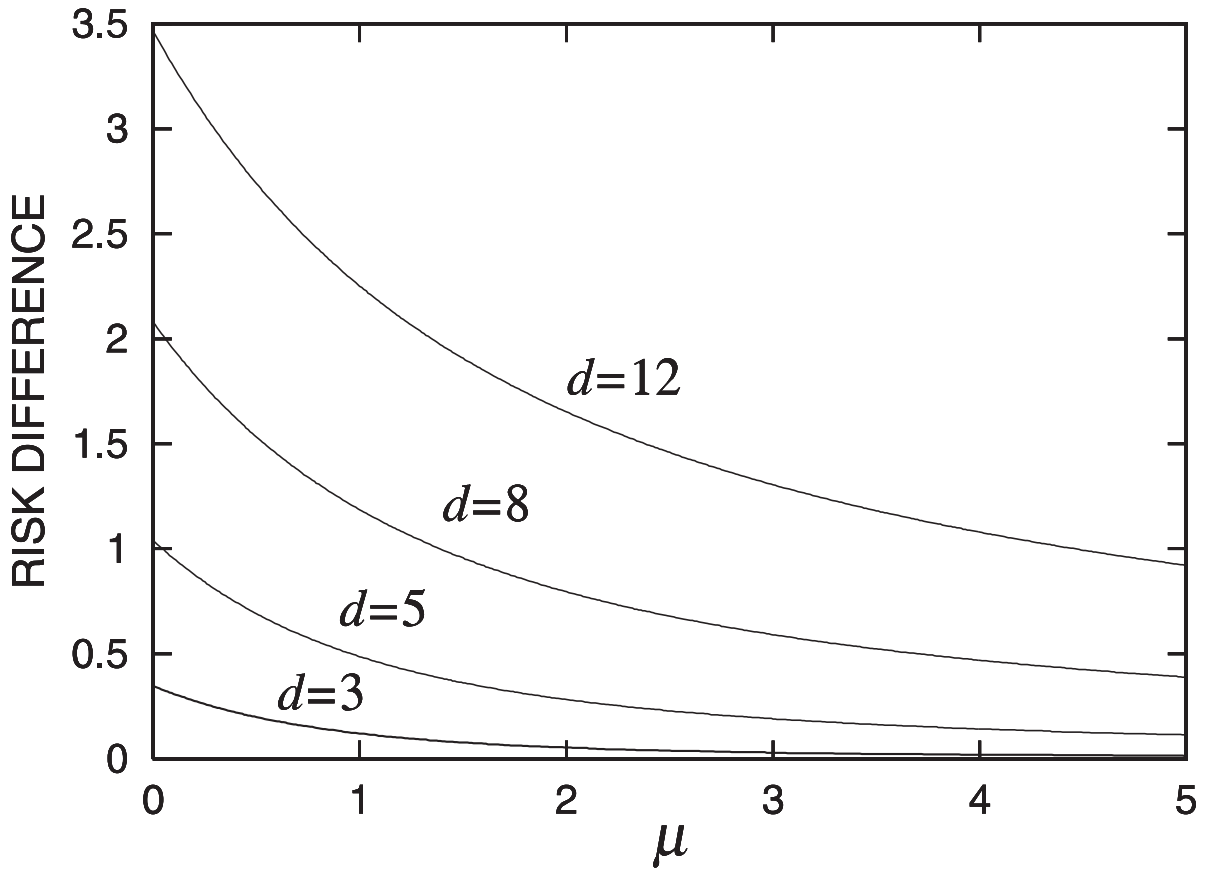

FIG. 1. *The difference between the expected divergences,* $E[D(p(y|\lambda), p_{\pi_J}(y|x))|\lambda] - E[D(p(y|\lambda), p_{\pi_S}(y|x))|\lambda]$, *which depends on* $\lambda$ *only through* $\mu = \lambda_1 + \lambda_2 + \cdots + \lambda_d$, *for* $d = 3, 5, 8, 12$.

When we adopt a prior distribution $\pi(\lambda)$, the plug-in distribution $p(y|\hat{\lambda}(x))$, where $\hat{\lambda}(x)$ is the generalized Bayes estimator based on $\pi(\lambda)$, is often used for prediction.

THEOREM 5. *The plug-in distribution* $p(y|\hat{\lambda}(x))$ *with the generalized Bayes estimator* $\hat{\lambda}(x)$ *based on* $\pi_S(\lambda)$ *is inadmissible under the Kullback–Leibler loss.*

The proof of Theorem 5 is given in the Appendix.

It can be shown that the plug-in distribution $p(y|\hat{\lambda})$ with the generalized Bayes estimator $\hat{\lambda}$ based on $\pi_S(\lambda)$ is admissible in the class of all plug-in distributions. However, it is inadmissible in the class of all predictive distributions. Therefore, it is not reasonable to restrict the class of predictive distributions to plug-in distributions.

**4. Some asymptotic properties and discussion.** In this section, we discuss the relation between the results in the previous sections and several previous studies on Bayesian theory from asymptotic viewpoints. Suppose that $x(1), x(2), \ldots, x(n), x(n+1), \ldots, x(n+m)$ are independent random variables from a true density $p(x|\theta)$ that belongs to a statistical model $\{p(x|\theta) \mid \theta \in \Theta\}$. The dimension of the parameter space $\Theta$ is $d$. Let $x^{(n)} = (x(1), x(2), \ldots, x(n))$ and $x_{(m)} = (x(n+1), x(n+2), \ldots, x(n+m))$. The objective is to construct a good Bayesian predictive distribution $p_\pi(x_{(m)}|x^{(n)})$ based on a prior $\pi$.

When $x(1), x(2), \ldots, x(n), x(n+1), \ldots, x(n+m)$ are independent sets of $d$ independent Poisson random variables with mean parameters $\lambda_1, \ldots, \lambda_d$, we consider a slight generalization of the problem introduced in Section 1. The objective is to predict $y$ that is a set of $d$ independent Poisson random variables with mean parameters $b\lambda_1, b\lambda_2, \ldots, b\lambda_d$, $b > 0$, by using an observation $x$ that is a set of $d$ independent Poisson random variables with mean parameters $a\lambda_1, a\lambda_2, \ldots, a\lambda_d$, $a > 0$. Here, $a$ and $b$ correspond to $n$ and $m$, respectively.

4.1. *Some asymptotics.* First, we consider the asymptotics where $a$ and $d$ are fixed and $b$ goes to infinity. In this subsection $d \geq 3$ is assumed. The asymptotics are closely related to the setup where $n = 0$ and $m \to \infty$, which has been studied in reference analysis, coding theory and prequential analysis as we will see in the next subsection.

When $\mu := \sum_i \lambda_i = 0$, we have

$$\begin{aligned} &E[D(p(y|\lambda), p_{\pi_{\mathrm{J}}}(y|x))|\lambda] - E[D(p(y|\lambda), p_{\pi_{\mathrm{S}}}(y|x))|\lambda] \\ &\quad = \left(\frac{d}{2} - 1\right) \log \frac{a+b}{a} > 0 \end{aligned}$$

from (7). Thus, the risk difference between the Bayesian predictive distribution based on the Jeffreys prior $\pi_{\mathrm{J}}(\lambda)$ and that based on the shrinkage prior $\pi_{\mathrm{S}}(\lambda)$ is of order $\log b$.

When $\mu \neq 0$, the risk difference converges to a positive constant when $a$ and $d$ are fixed and $b$ goes to infinity. By evaluating (6) using Stirling's formula $\log \Gamma(x) = \log(2\pi)^{1/2} + (x - 1/2)\log x - x + o(1)$, it can be easily verified that

$$\begin{aligned} &\lim_{b \to \infty} (E[D(p(y|\lambda), p_{\pi_{\mathrm{J}}}(y|x))|\lambda] - E[D(p(y|\lambda), p_{\pi_{\mathrm{S}}}(y|x))|\lambda]) \\ &\quad = E\left[\log \Gamma\left(\sum_i x_i + \frac{d}{2}\right) - \log \Gamma\left(\sum_i x_i + 1\right) - \frac{d-2}{2} \log a\mu \middle| \lambda\right] \\ &\quad > 0. \end{aligned}$$

Second, we consider the asymptotics where $b$ and $d$ are fixed and $a$ goes to infinity. There are many statistical applications where the objective is to construct a good predictive distribution for a future observation $x_{(m)}$ by using the observed data $x^{(n)}$ and $n$ is relatively large. An important example is one-step prediction, $b = m = 1$. Improper prior distributions are widely used to construct Bayesian predictive distributions. Asymptotic properties of predictive distributions for one-step prediction have been studied [Vidoni (1995), Komaki (1996), Hartigan (1998) and Komaki (2002b)]. When we consider the Poisson model, by a discussion similar to the previous studies, it can be shown that the loss function for a Bayesian predictive distribution

can be expanded as

$$
\begin{aligned}
&D(p(y|\theta), p_\pi(y|x)) \\
&\quad = \tfrac{1}{2}\sum_{i=1}^{d} g_{ii}(\theta)(\hat{\theta}^i_\pi - \theta^i)^2 + \text{terms independent of } \pi + O_p(a^{-2}),
\end{aligned} \tag{8}
$$

where $\theta^i = {\lambda_i}^{1/3}$, $\hat{\theta}^i_\pi = \{(\hat{\lambda}_\pi)_i\}^{1/3}$,

$$
(\hat{\lambda}_\pi)_i = \frac{x_i + 1/2}{a} + \left(\frac{x_i}{a^2}\right)\left\{\partial_{\lambda_i} \log\left(\frac{\pi(\lambda)}{\pi_{\mathrm{J}}(\lambda)}\right)\right\}\bigg|_{\lambda_i = x_i/a} + o_p(a^{-1})
$$

and $ag_{ii}(\theta) = 9a\theta^i$ is the Fisher information. The risk difference between the Bayesian predictive distribution based on the Jeffreys prior $\pi_{\mathrm{J}}(\lambda)$ and that based on the shrinkage prior $\pi_{\mathrm{S}}(\lambda)$ is of order $a^{-2}$ when $a$ goes to infinity [see Komaki (2002b) for details on the asymptotics of shrinkage predictive distributions]. Equation (8) gives an intuitive meaning for the Kullback–Leibler loss.

Third, we consider the asymptotics where $a$ and $b$ are fixed and $d$ goes to infinity. The data dimension $d$ becomes large in many fields of applied statistics such as spatial statistics, contingency table analysis and population data analysis.

It is easy to show the following result by evaluating (6) using Stirling's formula. If $\limsup_{d\to\infty}(\mu_d/d) < \infty$, where $\mu_d := \sum_{i=1}^d \lambda_i$, then

$$
\begin{aligned}
0 &< \liminf_{d\to\infty} \frac{1}{d}\{E[D(p(y|\lambda), p_{\pi_{\mathrm{J}}}(y|x))|\lambda] - E[D(p(y|\lambda), p_{\pi_{\mathrm{S}}}(y|x))|\lambda]\} \\
&\le \limsup_{d\to\infty} \frac{1}{d}\{E[D(p(y|\lambda), p_{\pi_{\mathrm{J}}}(y|x))|\lambda] - E[D(p(y|\lambda), p_{\pi_{\mathrm{S}}}(y|x))|\lambda]\} \\
&< \infty.
\end{aligned}
$$

For example, when $\lambda_i$, $i = 1, 2, 3, \ldots$, are generated independently from a distribution that has mean $\bar{\lambda}$, then $\lim_{d\to\infty}(\mu_d/d) = \bar{\lambda}$ almost surely and the risk difference is of order $d$ as $d$ goes to infinity.

4.2. *Relation to previous work.* In coding theory, the ideal code-word length of a Bayes code for a data string $x_{(m)} = (x(1), \ldots, x(m))$ based on a proper prior density $\pi(\theta)$ is given by

$$
\begin{aligned}
-\log p_\pi(x_{(m)}) &= -\log \int p(x(1), \ldots, x(m)|\theta)\pi(\theta)\, d\theta \\
&= -\sum_{l=0}^{m-1} \log \int p_\pi(x(l+1)|x_{(l)})\pi(\theta|x_{(l)})\, d\theta,
\end{aligned} \tag{9}
$$

where $n=0$. The average of the expected redundancy with respect to $\pi(\theta)$ is given by

$$\begin{aligned} & I_\pi(\theta; x_{(m)}) \\ &\quad = \int \pi(\theta) \int p(x_{(m)}|\theta) \biggl[ -\log\biggl\{ \int p(x_{(m)}|\theta')\pi(\theta')\, d\theta' \biggr\} \\ &\qquad\qquad + \log p(x_{(m)}|\theta) \biggr] \, dx_{(m)}\, d\theta \\ &\quad = \int \pi(\theta) \int p(x_{(m)}|\theta) \log \frac{p(x_{(m)}|\theta)\pi(\theta)}{\int p(x_{(m)}|\theta')\pi(\theta')\, d\theta' \, \pi(\theta)} \, dx_{(m)}\, d\theta, \end{aligned} \tag{10}$$

which is the mutual information between $\theta$ and $x_{(m)}$.

Bernardo (1979) introduced the notion of reference prior distributions and showed that the Jeffreys prior asymptotically maximizes the mutual information between $\theta$ and $x_{(m)} = (x(1), x(2), \dots, x(m))$ when $m \to \infty$ by using a heuristic discussion, although the mutual information cannot be properly defined if $\pi(\theta)$ is improper. Prequential analysis [Dawid (1984) and Skouras and Dawid (1999)] is also based on the logarithmic scoring rule used to give code lengths.

In the discussions above, there exist serious technical difficulties associated with infinite integrals when we consider improper prior distributions. If $\pi(\theta)$ is improper (9) cannot be regarded as an ideal code-word length of a Bayes code. A compact subset or a sequence of compact subsets of the original parameter space $\Theta$ has been considered to handle the difficulties in many previous studies. The heuristics are artificial but useful for treating the problems rigorously.

When $n=0$ and $m$ goes to infinity, under suitable regularity conditions, the mutual information between $x_{(m)}$ and $\theta$ is expanded as

$$\begin{aligned} I_\pi(\theta; x_{(m)}) &= \frac{d}{2} \log \frac{m}{2\pi e} + \int_K \pi(\theta) \log |g(\theta)|^{1/2} \, d\theta \\ &\quad - \int_K \pi(\theta) \log \pi(\theta)\, d\theta + o(1), \end{aligned} \tag{11}$$

where $K$ is a compact subset of the original parameter space $\Theta$ and $|g(\theta)|$ is the determinant of the Fisher information matrix [Ibragimov and Hasminskii (1973) and Clarke and Barron (1994)]. Thus (11) is maximized when $\pi(\theta) \propto |g(\theta)|^{1/2}$, which is the Jeffreys prior. The difference in $I_\pi(\theta; x_{(m)})$ due to the choice of a prior $\pi(\theta)$ is of order 1 when $m$ goes to infinity.

Here we consider the Poisson model and introduce an alternative method to deal with the difficulties associated with improper priors. Suppose that a transmitter $A$ and a receiver $B$ commonly observe a data sequence $x^{(n)} =$ $(x(1), x(2), \ldots, x(n))$. Only the transmitter $A$ can observe the subsequent data sequence $x_{(m)} = (x(n+1), \ldots, x(n+m))$. The transmitter $A$ sends $x_{(m)}$ to the receiver $B$ by using a Bayes code based on a prior $\pi(\lambda)$. Then the ideal code-word length of the Bayes code for $x_{(m)}$ can properly be defined by

$$-\log \int p(x_{(m)}|\lambda)\pi(\lambda|x^{(n)})\, d\lambda$$

if the posterior density $\pi(\theta|x^{(n)})$ is a proper density. The Bayes risk

$$\int \pi(\lambda) \sum_{x^{(n)}} p(x^{(n)}|\lambda) \sum_{x_{(m)}} p(x_{(m)}|\lambda) \log \frac{p(x_{(m)}|\lambda)}{p_\pi(x_{(m)}|x^{(n)})}\, d\lambda$$

coincides with the mutual information (10) when $n = 0$.

Now we consider the slightly generalized Poisson model. When $a$ is close to 0, the observation $x$ provides only a small amount of information and the situation is close to the setup that has been studied in reference analysis and Bayes coding theory, where the Jeffreys priors are recommended. However, Corollary 1 in Section 3 shows that the Bayesian predictive distribution based on the shrinkage prior $\pi_{\mathrm{S}}(\lambda)$ has better performance than that based on the Jeffreys prior $\pi_{\mathrm{J}}(\lambda)$ even in such a situation, since the risk function of the shrinkage prior is smaller than that of the Jeffreys prior for all $a > 0$ and $b > 0$ [see also Komaki (2002a) for related discussion for group models].

Note that our discussion is based on the original parameter space. It seems difficult to analyze the shrinkage phenomenon under the assumption that the real parameter value is in a compact subset of the original parameter space.

Finally, we note that predictive distributions based on the Jeffreys prior seem to become inadmissible under many loss functions other than the Kullback–Leibler loss. The admissibility of the predictive distribution based on the Jeffreys prior and shrinkage predictive distributions under other loss functions requires further study, although the Kullback–Leibler divergence is a natural loss function in several important streams in Bayesian theory.

## APPENDIX

Proofs of lemmas and Theorems 3 and 5.

Proof of Lemma 1. Let $\mu = \lambda_1 + \lambda_2 + \cdots + \lambda_d$ and $w_i = \lambda_i/\mu$, $i = 1, \ldots, d-1$. Since the relation

$$dw_1 \cdots dw_{d-1}\, d\mu = \left(\frac{1}{\mu}\right)^{d-1} d\lambda_1 \cdots d\lambda_d$$

holds, we have

$$\begin{aligned}
&\int \left(\sum_{i=1}^{d} \lambda_i\right)^{-\alpha} \lambda_1^{\beta_1 - 1} \lambda_2^{\beta_2 - 1} \cdots \lambda_d^{\beta_d - 1} \prod_{i=1}^{d} \{\exp(-a\lambda_i)(a\lambda_i)^{x_i}\}\, d\lambda_1\, d\lambda_2 \cdots d\lambda_d \\
&\quad = \int_0^\infty \mu^{-\alpha + \sum \beta_i + \sum x_i - 1} \exp(-a\mu) a^{\sum x_i} \\
&\qquad\quad \times w_1^{\beta_1 + x_1 - 1} w_2^{\beta_2 + x_2 - 1} \cdots w_{d-1}^{\beta_{d-1} + x_{d-1} - 1} \\
&\qquad\quad \times \left(1 - \sum_{i=1}^{d-1} w_i\right)^{\beta_d + x_d - 1} dw_1\, dw_2 \cdots dw_{d-1}\, d\mu \\
&\quad = a^{\alpha - \sum \beta_i} \frac{\Gamma(\sum_i x_i - \alpha + \sum_i \beta_i)}{\Gamma(\sum_i x_i + \sum_i \beta_i)} \prod_{i=1}^{d} \Gamma(x_i + \beta_i). \qquad \square
\end{aligned}$$

PROOF OF LEMMA 2. The derivative of $E[\log\Gamma(X+1+c) - \log\Gamma(X+1) - c\log\mu|\mu]$ satisfies the following inequality:

$$\begin{aligned}
&\frac{\partial}{\partial\mu} E[\log\Gamma(X+1+c) - \log\Gamma(X+1) - c\log\mu|\mu] \\
&\quad = -\sum_{k=0}^{\infty} \exp(-\mu)\frac{\mu^k}{k!} \log\Gamma(k+1+c) \\
&\qquad + \sum_{k=1}^{\infty} \exp(-\mu)\frac{\mu^{k-1}}{(k-1)!} \log\Gamma(k+1+c) \\
&\qquad + \sum_{k=0}^{\infty} \exp(-\mu)\frac{\mu^k}{k!} \log\Gamma(k+1) \\
&\qquad - \sum_{k=1}^{\infty} \exp(-\mu)\frac{\mu^{k-1}}{(k-1)!} \log\Gamma(k+1) - \frac{c}{\mu} \\
&\quad = \sum_{k=0}^{\infty} \exp(-\mu)\frac{\mu^k}{k!} \{\log(k+1+c) - \log(k+1)\} - \frac{c}{\mu} \\
&\quad \le \sum_{k=0}^{\infty} \exp(-\mu)\frac{\mu^k}{k!} \frac{c}{k+1} - \frac{c}{\mu} \\
&\quad = \frac{c}{\mu} \exp(-\mu) \left\{ \sum_{k=0}^{\infty} \frac{\mu^{k+1}}{(k+1)!} - \exp\mu \right\} \\
&\quad = -\frac{c}{\mu} \exp(-\mu) < 0.
\end{aligned}$$

We have thus proved the desired result. □

PROOF OF THEOREM 3. The admissibility is proved by using Blyth's method [Blyth (1951)]. For convenience, we put

$$\pi_{\alpha,\beta}(\lambda)\,d\lambda_1\,d\lambda_2\cdots d\lambda_d = \frac{\Gamma(\sum\beta_i)}{\prod_{i=1}^d \Gamma(\beta_i)}\frac{\lambda_1^{\beta_1-1}\lambda_2^{\beta_2-1}\cdots\lambda_d^{\beta_d-1}}{(\lambda_1+\lambda_2+\cdots+\lambda_d)^{\alpha}}\,d\lambda_1\,d\lambda_2\cdots d\lambda_d$$

in this proof. We use a sequence of priors $\{\pi^{[l]}_{\alpha,\beta}(\lambda) = \pi_{\alpha,\beta}(\lambda)\frac{1}{2}h_l^2(\mu)\}$, where $\{h_l\}$ is a sequence of functions defined by

$$h_l(\mu) = \begin{cases} 1, & \text{if } 0 \le \mu \le 1, \\ 1 - \dfrac{\log\mu}{\log l}, & \text{if } 1 < \mu \le l, \\ 0, & \text{if } l < \mu. \end{cases}$$

Function sequences of this kind are introduced by Brown and Hwang (1982) and have been used to prove the admissibility of various generalized Bayes estimators.

First, we see that the Bayesian predictive distribution

$$p_{\pi^{[l]}_{\alpha,\beta}}(y|x) = \frac{\int p(y|\lambda)p(x|\lambda)\pi^{[l]}_{\alpha,\beta}(\lambda)\,d\lambda}{\int p(x|\lambda)\pi^{[l]}_{\alpha,\beta}(\lambda)\,d\lambda}$$

based on the prior $\pi^{[l]}_{\alpha,\beta}(\lambda)$ minimizes the Bayes risk under $\pi^{[l]}_{\alpha,\beta}(\lambda)$ by using Aitchison's discussion [Aitchison (1975), page 549].

The Bayes risk of a predictive distribution $\hat{p}(y;x)$ is given by

$$\begin{aligned}
&\int \pi^{[l]}_{\alpha,\beta}(\lambda)\sum_x p(x|\lambda)\sum_y p(y|\lambda)\log\frac{p(y|\lambda)}{\hat{p}(y;x)}\,d\lambda \\
&\quad= \int \pi^{[l]}_{\alpha,\beta}(\lambda)\sum_x p(x|\lambda)\sum_y p(y|\lambda)\log p(y|\lambda)\,d\lambda \\
&\qquad- \int \pi^{[l]}_{\alpha,\beta}(\lambda)\sum_x p(x|\lambda)\sum_y p(y|\lambda)\log\hat{p}(y;x)\,d\lambda.
\end{aligned} \tag{12}$$

The first term on the right-hand side of (12) is finite and does not depend on $\hat{p}(y;x)$. The second term on the right-hand side of (12) is

$$\begin{aligned}
&-\int \pi^{[l]}_{\alpha,\beta}(\lambda)\sum_x p(x|\lambda)\sum_y p(y|\lambda)\log\hat{p}(y;x)\,d\lambda \\
&\quad= -\sum_x\sum_y p_{\pi^{[l]}_{\alpha,\beta}}(x,y)\log\hat{p}(y;x) \\
&\quad= -\sum_x p_{\pi^{[l]}_{\alpha,\beta}}(x)\sum_y p_{\pi^{[l]}_{\alpha,\beta}}(y|x)\log\hat{p}(y;x),
\end{aligned}$$

where

$$p_{\pi^{[l]}_{\alpha,\beta}}(x,y) = \int p(x|\lambda)p(y|\lambda)\pi^{[l]}_{\alpha,\beta}(\lambda)\,d\lambda$$

and

$$p_{\pi^{[l]}_{\alpha,\beta}}(x) = \int p(x|\lambda)\pi^{[l]}_{\alpha,\beta}(\lambda)\,d\lambda.$$

This is minimized when $\hat{p}(y;x) = p_{\pi^{[l]}_{\alpha,\beta}}(y|x)$. Thus, $p_{\pi^{[l]}_{\alpha,\beta}}(y|x)$ minimizes the Bayes risk (12).

Therefore, it suffices to show that

$$\begin{aligned}\int \pi^{[l]}_{\alpha,\beta}(\lambda)\Big\{&E[D(p(y|\lambda),p_{\pi_{\alpha,\beta}}(y|x))|\lambda] \\ &- E[D(p(y|\lambda),p_{\pi^{[l]}_{\alpha,\beta}}(y|x))|\lambda]\Big\}\,d\lambda \to 0 \qquad \text{as } l\to\infty\end{aligned} \tag{13}$$

to prove the admissibility of the Bayesian predictive distributions based on the priors in $\{\pi_{\alpha,\beta}(\lambda) : 0 < -\alpha + \sum_i \beta_i \le 1,\ \beta_i > 0,\ i = 1,2,\dots,d\}$ because the risks of the Bayesian predictive distributions with the priors in the proposed class are finite for all values of $\lambda$ and Theorem 2 holds.

Now we obtain a convenient expression for the integral

$$\int \pi^{[l]}_{\alpha,\beta}(\lambda)\{E[D(p(y|\lambda),p_{\pi_{\alpha,\beta}}(y|x))|\lambda] - E[D(p(y|\lambda),p_{\pi^{[l]}_{\alpha,\beta}}(y|x))|\lambda]\}\,d\lambda.$$

Let

$$\pi_{\alpha,\beta}(\mu) := \left(\frac{1}{\mu}\right)^{\alpha-\sum\beta_i+1}, \qquad \pi^{[l]}_{\alpha,\beta}(\mu) := \frac{1}{2}h_l^2(\mu)\left(\frac{1}{\mu}\right)^{\alpha-\sum\beta_i+1}$$

and

$$\begin{aligned}&\pi_{\alpha,\beta}(w_1,w_2,\dots,w_{d-1}) \\ &\quad := \pi^{[l]}_{\alpha,\beta}(w_1,w_2,\dots,w_{d-1}) \\ &\quad := \frac{\Gamma(\sum\beta_i)}{\prod_i\Gamma(\beta_i)} w_1^{\beta_1-1}w_2^{\beta_2-1}\cdots w_{d-1}^{\beta_{d-1}-1}\left(1-\sum_{i=1}^{d-1}w_i\right)^{\beta_d-1}.\end{aligned} \tag{14}$$

Then we have

$$\begin{aligned}&\pi_{\alpha,\beta}(\lambda)\,d\lambda_1\,d\lambda_2\cdots d\lambda_d \\ &\quad = \pi_{\alpha,\beta}(\mu,w_1,\dots,w_{d-1})\,d\mu\,dw_1\cdots dw_{d-1} \\ &\quad = \pi_{\alpha,\beta}(\mu)\pi_{\alpha,\beta}(w_1,\dots,w_{d-1})\,d\mu\,dw_1\cdots dw_{d-1}\end{aligned} \tag{15}$$

and

$$\begin{aligned}
&\pi^{[l]}_{\alpha,\beta}(\lambda)\,d\lambda_1\,d\lambda_2\cdots d\lambda_d\\
&\qquad=\pi^{[l]}_{\alpha,\beta}(\mu,w_1,\ldots,w_{d-1})\,d\mu\,dw_1\cdots dw_{d-1}\\
&\qquad=\pi^{[l]}_{\alpha,\beta}(\mu)\pi^{[l]}_{\alpha,\beta}(w_1,\ldots,w_{d-1})\,d\mu\,dw_1\cdots dw_{d-1}.
\end{aligned}\tag{16}$$

Let $\tilde{x}=x_1+x_2+\cdots+x_d$ and $\tilde{y}=y_1+y_2+\cdots+y_d$. If a prior $\pi(\mu,w_1,w_2,\ldots,w_{d-1})$ has the form $\pi(\mu,w_1,w_2,\ldots,w_{d-1})=\pi(\mu)\pi(w_1,w_2,\ldots,w_{d-1})$, then the relation

$$\begin{aligned}
&p_\pi(\mu,w_1,w_2,\ldots,w_{d-1}|x_1,x_2,\ldots,x_d)\\
&\quad=[p(\tilde{x},x_1,x_2,\ldots,x_{d-1}|\mu,w_1,w_2,\ldots,w_{d-1})\\
&\qquad\qquad\qquad\times\pi(\mu,w_1,w_2,\ldots,w_{d-1})]\\
&\qquad\times\biggl[\int p(\tilde{x},x_1,x_2,\ldots,x_{d-1}|\mu,w_1,w_2,\ldots,w_{d-1})\\
&\qquad\qquad\times\pi(\mu,w_1,w_2,\ldots,w_{d-1})\,d\mu\,dw_1\,dw_2\cdots dw_{d-1}\biggr]^{-1}\\
&\quad=[p(\tilde{x}|\mu)\pi(\mu)]\biggl[\int p(\tilde{x}|\mu)\pi(\mu)\,d\mu\biggr]^{-1}\\
&\qquad\times[p(x_1,x_2,\ldots,x_{d-1}|\tilde{x},w_1,w_2,\ldots,w_{d-1})\pi(w_1,w_2,\ldots,w_{d-1})]\\
&\qquad\times\biggl[\int p(x_1,x_2,\ldots,x_{d-1}|\tilde{x},w_1,w_2,\ldots,w_{d-1})\\
&\qquad\qquad\times\pi(w_1,w_2,\ldots,w_{d-1})\,dw_1\,dw_2\cdots dw_{d-1}\biggr]^{-1}
\end{aligned}\tag{17}$$

holds, because

$$\begin{aligned}
&p(x_1,x_2,\ldots,x_d|\lambda_1,\lambda_2,\ldots,\lambda_d)\\
&\qquad=p(\tilde{x},x_1,x_2,\ldots,x_{d-1}|\mu,w_1,w_2,\ldots,w_{d-1})\\
&\qquad=p(\tilde{x}|\mu,w_1,w_2,\ldots,w_{d-1})\\
&\qquad\quad\times p(x_1,x_2,\ldots,x_{d-1}|\tilde{x},\mu,w_1,w_2,\ldots,w_{d-1})\\
&\qquad=p(\tilde{x}|\mu)p(x_1,x_2,\ldots,x_{d-1}|\tilde{x},w_1,w_2,\ldots,w_{d-1}).
\end{aligned}$$

From the relations (14)–(17), and

$$\begin{aligned}
&p(y_1,y_2,\ldots,y_d|\lambda_1,\lambda_2,\ldots,\lambda_d)\\
&\qquad=p(\tilde{y}|\mu)p(y_1,y_2,\ldots,y_{d-1}|\tilde{y},w_1,w_2,\ldots,w_{d-1}),
\end{aligned}$$

it follows that the difference of the Kullback–Leibler losses for

$$p_{\pi_{\alpha,\beta}}(y_1, y_2, \ldots, y_d | x_1, x_2, \ldots, x_d)$$

and

$$p_{\pi^{[l]}_{\alpha,\beta}}(y_1, y_2, \ldots, y_d | x_1, x_2, \ldots, x_d)$$

is given by

$$\begin{aligned}
&D(p(y|\lambda), p_{\pi_{\alpha,\beta}}(y|x)) - D(p(y|\lambda), p_{\pi^{[l]}_{\alpha,\beta}}(y|x)) \\
&\quad = \sum_{y_1}\sum_{y_2}\cdots\sum_{y_d} p(y_1, y_2, \ldots, y_d | \lambda_1, \lambda_2, \ldots, \lambda_d) \\
&\qquad\qquad \times \log \frac{p(y_1, y_2, \ldots, y_d | \lambda_1, \lambda_2, \ldots, \lambda_d)}{p_{\pi_{\alpha,\beta}}(y_1, y_2, \ldots, y_d | x_1, x_2, \ldots, x_d)} \\
&\qquad - \sum_{y_1}\sum_{y_2}\cdots\sum_{y_d} p(y_1, y_2, \ldots, y_d | \lambda_1, \lambda_2, \ldots, \lambda_d) \\
&\qquad\qquad \times \log \frac{p(y_1, y_2, \ldots, y_d | \lambda_1, \lambda_2, \ldots, \lambda_d)}{p_{\pi^{[l]}_{\alpha,\beta}}(y_1, y_2, \ldots, y_d | x_1, x_2, \ldots, x_d)} \qquad (18) \\
&\quad = \sum_{y_1}\sum_{y_2}\cdots\sum_{y_d} p(\tilde{y}|\mu) p(y_1, y_2, \ldots, y_{d-1} | \tilde{y}, w_1, w_2, \ldots, w_{d-1}) \\
&\qquad\qquad \times \log \frac{p_{\pi^{[l]}_{\alpha,\beta}}(y_1, y_2, \ldots, y_d | x_1, x_2, \ldots, x_d)}{p_{\pi_{\alpha,\beta}}(y_1, y_2, \ldots, y_d | x_1, x_2, \ldots, x_d)} \\
&\quad = \sum_{\tilde{y}} p(\tilde{y}|\mu) \log \frac{p_{\pi^{[l]}_{\alpha,\beta}}(\tilde{y}|\tilde{x})}{p_{\pi_{\alpha,\beta}}(\tilde{y}|\tilde{x})},
\end{aligned}$$

where

$$p_{\pi_{\alpha,\beta}}(\tilde{y}|\tilde{x}) = \frac{\int_0^\infty p(\tilde{y}|\mu) p(\tilde{x}|\mu) \pi_{\alpha,\beta}(\mu)\, d\mu}{\int_0^\infty p(\tilde{x}|\mu) \pi_{\alpha,\beta}(\mu)\, d\mu}$$

and

$$p_{\pi^{[l]}_{\alpha,\beta}}(\tilde{y}|\tilde{x}) = \frac{\int_0^\infty p(\tilde{y}|\mu) p(\tilde{x}|\mu) \pi^{[l]}_{\alpha,\beta}(\mu)\, d\mu}{\int_0^\infty p(\tilde{x}|\mu) \pi^{[l]}_{\alpha,\beta}(\mu)\, d\mu}.$$

Since

$$\begin{aligned}
&\sum_{\tilde{x}} \exp(-a\mu) \frac{(a\mu)^{\tilde{x}}}{\tilde{x}!} \\
&\quad \times \sum_{\tilde{y}} \exp\{-(s+\tau)\mu\} \frac{\{(s+\tau)\mu\}^{\tilde{y}}}{\tilde{y}!}
\end{aligned}$$

$$\times \log\left[\exp\{-(s+\tau)\mu\}\frac{\{(s+\tau)\mu\}^{\tilde y}}{\tilde y!}\right.$$
$$\left.\times\left(\int_0^\infty \exp\{-(s+\tau)\mu\}\frac{\{(s+\tau)\mu\}^{\tilde y}}{\tilde y!}\pi(\mu|\tilde x)\,d\mu\right)^{-1}\right]$$
$$=\sum_{\tilde x}\exp(-a\mu)\frac{(a\mu)^{\tilde x}}{\tilde x!}$$
$$\times\sum_v\sum_w \exp(-s\mu)\frac{(s\mu)^v}{v!}\exp(-\tau\mu)\frac{(\tau\mu)^w}{w!}$$
$$\times\left\{\log\frac{\exp(-s\mu)\mu^v}{\int_0^\infty \exp(-s\mu)\mu^v\pi(\mu|\tilde x)\,d\mu}\right.$$
$$+\log\left[(\exp(-\tau\mu)\mu^w)\right.$$
$$\left.\left.\times\left(\frac{\int_0^\infty \exp(-s\mu)\mu^v\pi(\mu|\tilde x)\,d\mu}{\int_0^\infty \exp(-s\mu)\mu^v\exp(-\tau\mu)\mu^w\pi(\mu|\tilde x)\,d\mu}\right)\right]\right\},$$

we have

$$\frac{d}{ds}\sum_{\tilde x}\exp(-a\mu)\frac{(a\mu)^{\tilde x}}{\tilde x!}$$
$$\times\sum_{\tilde y}\exp(-s\mu)\frac{(s\mu)^{\tilde y}}{\tilde y!}\log\frac{\exp(-s\mu)(s\mu)^{\tilde y}/(\tilde y!)}{\int_0^\infty \exp(-s\mu)(s\mu)^{\tilde y}/(\tilde y!)\pi(\mu|\tilde x)\,d\mu}$$
$$=\lim_{\tau\to 0}\frac{1}{\tau}\left[\sum_{\tilde x}\sum_v\sum_w \exp(-a\mu)\frac{(a\mu)^{\tilde x}}{\tilde x!}\exp(-s\mu)\frac{(s\mu)^v}{v!}\exp(-\tau\mu)\frac{(\tau\mu)^w}{w!}\right.$$
$$\left.\times\log\left[\frac{\exp(-\tau\mu)\mu^w\int_0^\infty \exp(-s\mu)\mu^v\pi(\mu|\tilde x)\,d\mu}{\int_0^\infty \exp(-s\mu)\mu^v\exp(-\tau\mu)\mu^w\pi(\mu|\tilde x)\,d\mu}\right]\right]$$
$$=\sum_z \exp(-t\mu)\frac{(t\mu)^z}{z!}\left(\hat\mu-\mu-\mu\log\frac{\hat\mu}{\mu}\right),$$

where $t := a+s$, $z := \tilde x + v$ and

$$\hat\mu := \frac{\int_0^\infty \mu\exp(-s\mu)(s\mu)^v/(v!)\pi(\mu|\tilde x)\,d\mu}{\int_0^\infty \exp(-s\mu)(s\mu)^v/(v!)\pi(\mu|\tilde x)\,d\mu}$$
$$=\frac{\int_0^\infty \mu\exp(-t\mu)(t\mu)^z/(z!)\pi(\mu)\,d\mu}{\int_0^\infty \exp(-t\mu)(t\mu)^z/(z!)\pi(\mu)\,d\mu}.$$

We put $c = \alpha - \sum_i \beta_i + 1$ $(0 \le c < 1)$ and $g_l(\mu) = (1/2)h_l^2(\mu)$. Then the expected divergence from $p(\tilde{y}|\mu)$ to $p_{\pi^{[l]}_{\alpha,\beta}}(\tilde{y}|\tilde{x})$ is expressed by

$$
\begin{aligned}
&E[D(p(\tilde{y}|\mu), p_{\pi^{[l]}_{\alpha,\beta}}(\tilde{y}|\tilde{x}))|\mu] \\
&\quad = \sum_{\tilde{x}} \exp(-a\mu) \frac{(a\mu)^{\tilde{x}}}{\tilde{x}!} \\
&\qquad \times \sum_{\tilde{y}} \exp(-b\mu) \frac{(b\mu)^{\tilde{y}}}{\tilde{y}!} \log \frac{\exp(-b\mu)(b\mu)^{\tilde{y}}/(\tilde{y}!)}{\int \exp(-b\bar{\mu})(b\bar{\mu})^{\tilde{y}}/(\tilde{y}!) p_{\pi^{[l]}_{\alpha,\beta}}(\bar{\mu}|\tilde{x})\, d\bar{\mu}} \\
&\quad = \int_a^{a+b} \sum_z \exp(-t\mu) \frac{(t\mu)^z}{z!} \left( \hat{\mu}_{l,t} - \mu - \mu \log \frac{\hat{\mu}_{l,t}}{\mu} \right) dt,
\end{aligned} \tag{19}
$$

where

$$
\begin{aligned}
\hat{\mu}_{l,t} &= \frac{\int_0^\infty \mu \exp(-t\mu)(t\mu)^z/(z!)\mu^{-c} g_l(\mu)\, d\mu}{\int_0^\infty \exp(-t\mu)(t\mu)^z/(z!)\mu^{-c} g_l(\mu)\, d\mu} \\
&= \frac{\int_0^\infty \exp(-t\mu)(t\mu)^{z+1-c} g_l(\mu)\, d\mu}{t \int_0^\infty \exp(-t\mu)(t\mu)^{z-c} g_l(\mu)\, d\mu} \\
&= \Big[ -t^{-1} \exp(-t\mu)(t\mu)^{z+1-c} g_l(\mu) \big|_0^\infty \\
&\qquad + \int_0^\infty t^{-1} \exp(-t\mu) \{ (z+1-c)(t\mu)^{z-c} t g_l(\mu) \\
&\qquad\qquad + (t\mu)^{z+1-c} g_l'(\mu) \}\, d\mu \Big] \\
&\qquad \times \left[ t \int_0^\infty \exp(-t\mu)(t\mu)^{z-c} g_l(\mu)\, d\mu \right]^{-1} \\
&= \frac{z+1-c}{t} + \frac{\int_0^\infty \exp(-t\mu)(t\mu)^{z+1-c} g_l'(\mu)\, d\mu}{t^2 \int_0^\infty \exp(-t\mu)(t\mu)^{z-c} g_l(\mu)\, d\mu}.
\end{aligned} \tag{20}
$$

In the same way, the expected divergence from $p(\tilde{y}|\mu)$ to

$$
p_{\pi_{\alpha,\beta}}(\tilde{y}|\tilde{x}) = \left( \frac{a}{a+b} \right)^{\tilde{x}+1-c} \left( \frac{b}{a+b} \right)^{\tilde{y}} \frac{\Gamma(\tilde{x}+\tilde{y}-c+1)}{\Gamma(\tilde{x}-c+1)\tilde{y}!}
$$

is expressed by

$$
\begin{aligned}
&E[D(p(\tilde{y}|\mu), p_{\pi_{\alpha,\beta}}(\tilde{y}|\tilde{x}))|\mu] \\
&\quad = \sum_{\tilde{x}} \exp(-a\mu) \frac{(a\mu)^{\tilde{x}}}{\tilde{x}!} \sum_{\tilde{y}} \exp(-b\mu) \frac{(b\mu)^{\tilde{y}}}{\tilde{y}!} \log \frac{\exp(-b\mu)(b\mu)^{\tilde{y}}/(\tilde{y}!)}{p_{\pi_{\alpha,\beta}}(\tilde{y}|\tilde{x})}
\end{aligned} \tag{21}
$$

$$= \int_a^{a+b} \sum_z \exp(-t\mu) \frac{(t\mu)^z}{z!} \left( \frac{z+1-c}{t} - \mu - \mu \log \frac{z+1-c}{t\mu} \right) dt.$$

From (18), (19) and (21), we obtain the following expression for the integral in (13):

$$\begin{aligned}
&\int \pi^{[l]}_{\alpha,\beta}(\lambda) \Big\{ E[D(p(y|\lambda), p_{\pi_{\alpha,\beta}}(y|x))|\lambda] \\
&\qquad - E[D(p(y|\lambda), p_{\pi^{[l]}_{\alpha,\beta}}(y|x))|\lambda] \Big\}\, d\lambda \\
&\quad = \int_0^\infty \pi^{[l]}_{\alpha,\beta}(\mu) \Big\{ E[D(p(\tilde{y}|\mu), p_{\pi_{\alpha,\beta}}(\tilde{y}|\tilde{x}))|\mu] \\
&\qquad\qquad - E[D(p(\tilde{y}|\mu), p_{\pi^{[l]}_{\alpha,\beta}}(\tilde{y}|\tilde{x}))|\mu] \Big\}\, d\mu \\
&\quad = \int_a^{a+b} \Big\{ \int_0^\infty g_l(\mu)\mu^{-c} \\
&\qquad \times \sum_{z=0}^{\infty} \exp(-t\mu) \frac{(t\mu)^z}{z!} \\
&\qquad\qquad \times \left( \frac{z+1-c}{t} - \hat{\mu}_{l,t} - \mu \log \frac{z+1-c}{t\hat{\mu}_{l,t}} \right) d\mu \Big\}\, dt.
\end{aligned} \tag{22}$$

We show (13) by evaluating (22) using the following inequalities, (23) and (24), similar to the inequalities used by Ghosh and Yang (1988) to prove the admissibility of a class of linear estimators of Poisson means of the form $\hat{\lambda}_i = c_i x_i + b_i$ under the Kullback–Leibler loss.

We have

$$\begin{aligned}
&\int_0^\infty g_l(\mu)\mu^{-c} \sum_{z=0}^{\infty} \exp(-t\mu) \frac{(t\mu)^z}{z!} \\
&\qquad \times \left( \frac{z+1-c}{t} - \hat{\mu}_{l,t} - \mu \log \frac{z+1-c}{t\hat{\mu}_{l,t}} \right) d\mu \\
&\quad \le \int_0^\infty g_l(\mu)\mu^{-c} \sum_{z=0}^{\infty} \exp(-t\mu) \frac{(t\mu)^z}{z!} \\
&\qquad\qquad \times \left\{ \frac{z+1-c}{t} - \hat{\mu}_{l,t} + \mu \frac{t\hat{\mu}_{l,t} - (z+1-c)}{z+1-c} \right\} d\mu \\
&\quad = \sum_{z=0}^{\infty} \frac{t^{c-1}}{z!} \Big\{ \int_0^\infty \exp(-t\mu) g_l(\mu) \frac{(t\mu)^{z+1-c}}{z+1-c}\, d\mu \\
&\qquad - \int_0^\infty \exp(-t\mu) g_l(\mu) (t\mu)^{z-c}\, d\mu \Big\} \{ t\hat{\mu}_{l,t} - (z+1-c) \}
\end{aligned} \tag{23}$$

$$
\begin{aligned}
&= \sum_{z=0}^{\infty} \frac{t^{c-1}}{z!} \biggl[ -t^{-1} \exp(-t\mu) g_l(\mu) \frac{(t\mu)^{z+1-c}}{z+1-c} \biggr|_0^{\infty} \\
&\qquad\qquad + \int_0^{\infty} t^{-1} \exp(-t\mu) \biggl\{ g_l'(\mu) \frac{(t\mu)^{z+1-c}}{z+1-c} + g_l(\mu) t (t\mu)^{z-c} \biggr\} \, d\mu \\
&\qquad\qquad\qquad\qquad - \int_0^{\infty} \exp(-t\mu) g_l(\mu) (t\mu)^{z-c} \, d\mu \biggr] \\
&\quad \times \{ t\hat{\mu}_{l,t} - (z+1-c) \} \\
&= \sum_{z=0}^{\infty} \frac{t^{c-2}}{z!} \biggl\{ \int_0^{\infty} \exp(-t\mu) g_l'(\mu) \frac{(t\mu)^{z+1-c}}{z+1-c} \, d\mu \biggr\} \{ t\hat{\mu}_{l,t} - (z+1-c) \}.
\end{aligned}
$$

By using (20) and the inequality

$$
\frac{z+1}{z+1-c} = \frac{1}{1-c},
$$

where $0 \leq c < 1$, we have

$$
\begin{aligned}
&\int_0^{\infty} g_l(\mu) \mu^{-c} \sum_{z=0}^{\infty} \exp(-t\mu) \frac{(t\mu)^z}{z!} \biggl( \frac{z+1-c}{t} - \hat{\mu}_{l,t} - \mu \log \frac{z+1-c}{t\hat{\mu}_{l,t}} \biggr) \, d\mu \\
&\leq \frac{t^{c-3}}{1-c} \sum_{z=0}^{\infty} \frac{1}{(z+1)!} \frac{\{ \int_0^{\infty} \exp(-t\mu)(t\mu)^{z+1-c} g_l'(\mu) \, d\mu \}^2}{\int \exp(-t\bar{\mu})(t\bar{\mu})^{z-c} g_l(\bar{\mu}) \, d\bar{\mu}} \qquad (24) \\
&= \frac{t^{c-3}}{1-c} \sum_{z=0}^{\infty} \frac{2}{(z+1)!} \frac{\{ \int_0^{\infty} \exp(-t\mu)(t\mu)^{z-c} (t\mu h_l'(\mu)) h_l(\mu) \, d\mu \}^2}{\int \exp(-t\bar{\mu})(t\bar{\mu})^{z-c} h_l^2(\bar{\mu}) \, d\bar{\mu}} \\
&\leq \frac{t^{c-3}}{1-c} \sum_{z=0}^{\infty} \frac{2}{(z+1)!} \frac{\int_0^{\infty} \exp(-t\mu)(t\mu)^{z-c} (t\mu h_l'(\mu))^2 \, d\mu}{\int \exp(-t\bar{\mu})(t\bar{\mu})^{z-c} h_l^2(\bar{\mu}) \, d\bar{\mu}} \\
&\qquad\qquad \times \int_0^{\infty} \exp(-t\mu)(t\mu)^{z-c} h_l^2(\mu) \, d\mu \\
&= \frac{t^{c-3}}{1-c} \sum_{z=0}^{\infty} \frac{2}{(z+1)!} \int_0^{\infty} \exp(-t\mu)(t\mu)^{z+2-c} (h_l'(\mu))^2 \, d\mu \\
&= \frac{2t^{c-3}}{1-c} \int_0^{\infty} \{1 - \exp(-t\mu)\} (t\mu)^{1-c} (h_l'(\mu))^2 \, d\mu \\
&\leq \frac{2}{(1-c)t^2} \int_0^{\infty} \mu^{1-c} (h_l'(\mu))^2 \, d\mu.
\end{aligned}
$$

The derivative of $h_l(\mu)$ is given by

$$
h_l'(\mu) = \begin{cases} 0, & \text{if } 0 < \mu < 1, \\ -\dfrac{1}{\mu \log l}, & \text{if } 1 < \mu < l, \\ 0, & \text{if } l < \mu. \end{cases} \qquad (25)
$$

From (22), (24) and (25), we have

$$\begin{aligned}&\int \pi^{[l]}_{\alpha,\beta}(\lambda)\{E[D(p(y|\mu), p_{\pi_{\alpha,\beta}}(y|x)) - D(p(y|\mu), p_{\pi^{[l]}_{\alpha,\beta}}(y|x))|\lambda]\}\,d\lambda\\
&\quad\le \int_a^{a+b}\left\{\frac{2}{(1-c)t^2}\int_0^\infty \mu^{1-c}(h_l'(\mu))^2\,d\mu\right\}dt\\
&\quad= \int_a^{a+b}\left\{\frac{2}{(1-c)t^2}\int_1^l \frac{1}{\mu^{1+c}(\log l)^2}\,d\mu\right\}dt\\
&\quad\le \int_a^{a+b}\frac{2}{(1-c)t^2}\frac{1}{\log l}\,dt = \frac{2}{(1-c)\log l}\left(\frac{1}{a}-\frac{1}{a+b}\right)\to 0\\
&\hspace{25em}\text{as } l\to\infty.\end{aligned}$$

We have thus proved the theorem. □

Proof of Theorem 5. From Lemma 1, the generalized Bayes estimator of $\lambda$ with respect to $\pi_{\mathrm{S}}(\lambda)$ is given by

$$\begin{aligned}\hat\lambda_i &= \frac{\int \lambda_i\pi_{\mathrm{S}}(\lambda)\prod_j\{\exp(-a\lambda_i)(a\lambda_i)^{x_i}/(x_i!)\}\,d\lambda}{\int \pi_{\mathrm{S}}(\lambda)\prod_j\{\exp(-a\lambda_i)(a\lambda_i)^{x_i}/(x_i!)\}\,d\lambda}\\
&= \frac{a^{-2}[\Gamma(\sum_k x_k+2)/\Gamma(\sum_l x_l + d/2+1)]\prod_{j\neq i}\Gamma(x_j+1/2)\Gamma(x_i+3/2)}{a^{-1}[\Gamma(\sum_k x_k+1)/\Gamma(\sum_l x_l+d/2)]\prod_j\Gamma(x_j+1/2)}\\
&= \frac{1}{a}\frac{\sum_j x_j+1}{\sum_k x_k+d/2}\left(x_i+\frac{1}{2}\right).\end{aligned}$$

The plug-in distribution with $\hat\lambda(x)$ is given by

$$p(y|\hat\lambda) = p(y|\hat\mu,\hat w) = p(\tilde y|\hat\mu)p(y|\tilde y,\hat w),$$

where

$$\hat\mu = \frac{1}{a}\left(\sum_i x_i + 1\right),\qquad \hat w_i = \frac{x_i+1/2}{\sum_j x_j + d/2},$$

$$p(\tilde y|\hat\mu) = \exp(-\hat\mu)\frac{\hat\mu^{\tilde y}}{\tilde y!}\quad\text{and}\quad p(y|\tilde y,\hat w) = \binom{\tilde y}{y_1 y_2\cdots y_d}\hat w_1^{y_1}\hat w_2^{y_2}\cdots\hat w_d^{y_d}.$$

We show that the predictive distribution $p_{\pi_{\mathrm{S}}}(\tilde y|\tilde x)p(y|\tilde y,\hat w)$ dominates the plug-in distribution $p(y|\hat\lambda)$. The difference between the risk of the plug-in distribution $p(y|\hat\lambda)$ and that of $p_{\pi_{\mathrm{S}}}(\tilde y|\tilde x)p(y|\tilde y,\hat w)$ is given by

$$\begin{aligned}&E[D(p(y|\lambda), p(y|\hat\lambda))|\lambda] - E[D(p(y|\lambda), p_{\pi_{\mathrm{S}}}(\tilde y|\tilde x)p(y|\tilde y,\hat w))|\lambda]\\
&\quad= E[D(p(\tilde y|\mu), p(\tilde y|\hat\mu))|\lambda] - E[D(p(y|\mu), p_{\pi_{\mathrm{S}}}(\tilde y|\tilde x))|\lambda].\end{aligned}$$

From (21), the expected Kullback–Leibler divergence from $p(\tilde{y}|\mu)$ to $p_{\pi_{\mathrm{S}}}(\tilde{y}|\tilde{x})$ is

$$
\begin{aligned}
&E[D(p(\tilde{y}|\mu), p_{\pi_{\mathrm{S}}}(\tilde{y}|\tilde{x}))|\lambda] \\
&\quad = \int_a^{a+b} \sum_{z=0}^{\infty} \exp(-t\mu) \frac{(t\mu)^z}{z!} \left( \frac{z+1}{t} - \mu - \mu \log \frac{z+1}{t\mu} \right) dt.
\end{aligned} \tag{26}
$$

Since the Kullback–Leibler divergence from $p(\tilde{y}|\mu)$ to $p(\tilde{y}|\hat{\mu})$ is given by

$$
D(p(\tilde{y}|\mu), p(\tilde{y}|\hat{\mu})) = b\left( \hat{\mu} - \mu - \mu \log \frac{\hat{\mu}}{\mu} \right),
$$

we have

$$
\begin{aligned}
&E[D(p(\tilde{y}|\mu), p(\tilde{y}|\hat{\mu}))|\lambda] \\
&\quad = b \sum_{\tilde{x}=0}^{\infty} \left\{ \exp(-a\mu) \frac{(a\mu)^{\tilde{x}}}{\tilde{x}!} \left( \frac{\tilde{x}+1}{a} - \mu - \mu \log \frac{\tilde{x}+1}{a\mu} \right) \right\}.
\end{aligned} \tag{27}
$$

Note that the integrand of (26) coincides with (27) multiplied by $1/b$ when $t = a$. Hence, to prove the theorem, it suffices to show that the integrand of (26) is a decreasing function of $t$ for all values of $\mu$.

The derivative of the integrand of (26) is given by

$$
\begin{aligned}
&\frac{d}{dt} \sum_{z=0}^{\infty} \exp(-t\mu) \frac{(t\mu)^z}{z!} \left( \frac{z+1}{t} - \mu - \mu \log \frac{z+1}{t\mu} \right) \\
&\quad = \mu^2 \left\{ -\frac{1}{t^2\mu^2} + \frac{1}{t\mu} + \sum_{z=0}^{\infty} \exp(-t\mu) \frac{(t\mu)^z}{z!} \log(z+1) \right. \\
&\qquad \left. - \sum_{z=1}^{\infty} \exp(-t\mu) \frac{(t\mu)^{z-1}}{(z-1)!} \log(z+1) \right\}.
\end{aligned}
$$

Since

$$
\begin{aligned}
-\frac{1}{t^2\mu^2} &= -\frac{1}{t^2\mu^2} \sum_{z=0}^{\infty} \exp(-t\mu) \frac{(t\mu)^z}{z!} \\
&= -\frac{1}{t^2\mu^2} \left\{ \exp(-t\mu) + t\mu \exp(-t\mu) + \sum_{z=0}^{\infty} \exp(-t\mu) \frac{(t\mu)^{z+2}}{(z+2)!} \right\} \\
&= -\frac{1}{t^2\mu^2} \{ \exp(-t\mu) + t\mu \exp(-t\mu) \} \\
&\quad - \sum_{z=0}^{\infty} \exp(-t\mu) \frac{(t\mu)^z}{z!} \left( \frac{1}{z+1} - \frac{1}{z+2} \right),
\end{aligned}
$$

$$\begin{aligned}\frac{1}{t\mu} &= \frac{1}{t\mu}\sum_{z=0}^{\infty}\exp(-t\mu)\frac{(t\mu)^z}{z!}\\ &= \frac{1}{t\mu}\exp(-t\mu) + \sum_{z=0}^{\infty}\exp(-t\mu)\frac{(t\mu)^z}{z!}\frac{1}{z+1}\end{aligned}$$

and

$$\begin{aligned}&\sum_{z=0}^{\infty}\exp(-t\mu)\frac{(t\mu)^z}{z!}\log(z+1) - \sum_{z=1}^{\infty}\exp(-t\mu)\frac{(t\mu)^{z-1}}{(z-1)!}\log(z+1)\\ &\quad = -\sum_{z=0}^{\infty}\exp(-t\mu)\frac{(t\mu)^z}{z!}\{\log(z+2)-\log(z+1)\}\\ &\quad \le -\sum_{z=0}^{\infty}\exp(-t\mu)\frac{(t\mu)^z}{z!}\frac{1}{z+2},\end{aligned}$$

we have

$$\frac{d}{dt}\sum_{z=0}^{\infty}\exp(-t\mu)\frac{(t\mu)^z}{z!}\left(\frac{z+1}{t} - \mu - \mu\log\frac{z+1}{t\mu}\right) \le -\frac{1}{t^2}\exp(-t\mu) < 0.$$

Thus, the integrand of (26) is a strictly decreasing function of $t$. Therefore, (26) is smaller than (27) for all values of $\lambda$. □

**Acknowledgment.** I am grateful to the referee for helpful and constructive comments.

## REFERENCES

Aitchison, J. (1975). Goodness of prediction fit. *Biometrika* **62** 547–554. MR391353

Aitchison, J. and Dunsmore, I. R. (1975). *Statistical Prediction Analysis.* Cambridge Univ. Press. MR408097

Akaike, H. (1978). A new look at the Bayes procedure. *Biometrika* **65** 53–59. MR501450

Berger, J. O. and Bernardo, J. M. (1989). Estimating a product of means: Bayesian analysis with reference priors. *J. Amer. Statist. Assoc.* **84** 200–207. MR999679

Bernardo, J. M. (1979). Reference posterior distributions for Bayesian inference (with discussion). *J. Roy. Statist. Soc. Ser. B* **41** 113–147. MR547240

Blyth, C. R. (1951). On minimax statistical decision procedures and their admissibility. *Ann. Math. Statist.* **22** 22–42. MR39966

Brown, L. D. and Hwang, J. T. (1982). A unified admissibility proof. In *Statistical Decision Theory and Related Topics III* (S. S. Gupta and J. O. Berger, eds.) **1** 205–230. Academic Press, New York. MR705290

Clarke, B. S. and Barron, A. R. (1994). Jeffreys' prior is asymptotically least favorable under entropy risk. *J. Statist. Plann. Inference* **41** 37–60. MR1292146

Clevenson, M. L. and Zidek, J. V. (1975). Simultaneous estimation of the means of independent Poisson laws. *J. Amer. Statist. Assoc.* **70** 698–705. MR394962

Dawid, A. P. (1984). Statistical theory. The prequential approach (with discussion). *J. Roy. Statist. Soc. Ser. A* **147** 278–292. MR763811

GEISSER, S. (1993). *Predictive Inference*: *An Introduction.* Chapman and Hall, New York. MR1252174

GHOSH, M. and YANG, M.-C. (1988). Simultaneous estimation of Poisson means under entropy loss. *Ann. Statist.* **16** 278–291. MR924871

HARTIGAN, J. A. (1965). The asymptotically unbiased prior distribution. *Ann. Math. Statist.* **36** 1137–1152. MR176539

HARTIGAN, J. A. (1998). The maximum likelihood prior. *Ann. Statist.* **26** 2083–2103. MR1700222

HAUSSLER, D. and OPPER, M. (1997). Mutual information, metric entropy and cumulative relative entropy risk. *Ann. Statist.* **25** 2451–2492. MR1604481

IBRAGIMOV, I. A. and HASMINSKII, R. Z. (1973). On the information contained in a sample about a parameter. In *Second International Symposium on Information Theory* (B. N. Petrov and F. Csáki, eds.) 295–309. Akadémiai Kiado, Budapest. MR356948

JAMES, W. and STEIN, C. (1961). Estimation with quadratic loss. *Proc. Fourth Berkeley Symp. Math. Statist. Probab.* **1** 361–380. Univ. California Press, Berkeley. MR133191

JEFFREYS, H. (1961). *Theory of Probability*, 3rd ed. Oxford Univ. Press. MR187257

KOMAKI, F. (1996). On asymptotic properties of predictive distributions. *Biometrika* **83** 299–313. MR1439785

KOMAKI, F. (2001). A shrinkage predictive distribution for multivariate normal observables. *Biometrika* **88** 859–864. MR1859415

KOMAKI, F. (2002a). Bayesian predictive distribution with right invariant priors. *Calcutta Statist. Assoc. Bull.* **52** 171–179. MR1971346

KOMAKI, F. (2002b). Shrinkage priors for Bayesian prediction. Unpublished manuscript.

MURRAY, G. D. (1977). A note on the estimation of probability density functions. *Biometrika* **64** 150–152. MR448690

NG, V. M. (1980). On the estimation of parametric density functions. *Biometrika* **67** 505–506. MR581751

SKOURAS, K. and DAWID, A. P. (1999). On efficient probability forecasting systems. *Biometrika* **86** 765–784. MR1741976

STEIN, C. (1956). Inadmissibility of the usual estimator for the mean of a multivariate distribution. *Proc. Third Berkeley Symp. Math. Statist. Probab.* **1** 197–206. Univ. California Press, Berkeley. MR84922

STEIN, C. (1974). Estimation of the mean of a multivariate normal distribution. In *Proc. Prague Symposium on Asymptotic Statistics* (J. Hájek, ed.) **2** 345–381. Univ. Karlova, Prague. MR381062

VIDONI, P. (1995). A simple predictive density based on the $p^*$-formula. *Biometrika* **82** 855–863. MR1380820

DEPARTMENT OF MATHEMATICAL INFORMATICS
GRADUATE SCHOOL OF INFORMATION
SCIENCE AND TECHNOLOGY
UNIVERSITY OF TOKYO
7-3-1 HONGO, BUNKYO-KU
TOKYO 113-8656
JAPAN
E-MAIL: komaki@mist.i.u-tokyo.ac.jp